\documentclass[11pt]{article}
\usepackage{amsmath}
\usepackage{amssymb}
\usepackage{amsfonts}
\usepackage{mathrsfs}
\usepackage{shortvrb,psfrag}
\usepackage{enumerate}
\usepackage{color}
\usepackage[dvips]{graphicx}

\let\pa=\partial

\let\f=\frac

\let\Om=\Omega

\def\R{\Bbb R}

\def\no{\noindent}
\def\na{\nabla}

\def\eqdef{\buildrel\hbox{\footnotesize def}\over =}
\def\endproof{\hphantom{MM}\hfill\llap{$\square$}\goodbreak}

\newcommand{\beq}{\begin{equation}}
\newcommand{\eeq}{\end{equation}}
\newcommand{\ben}{\begin{eqnarray}}
\newcommand{\een}{\end{eqnarray}}
\newcommand{\beno}{\begin{eqnarray*}}
\newcommand{\eeno}{\end{eqnarray*}}

\newtheorem{Theorem}{Theorem}[section]
\newtheorem{Definition}[Theorem]{Definition}
\newtheorem{Proposition}[Theorem]{Proposition}
\newtheorem{Lemma}[Theorem]{Lemma}

\newtheorem{Remark}[Theorem]{Remark}

\begin{document}
\title{On the interior regularity criteria of the 3-D Navier-Stokes equations involving two velocity components}

\author{Wendong Wang$^\dag$\, Liqun Zhang$^\ddag$\,and\, Zhifei Zhang$^\sharp$\\[2mm]
{\small $^\dag$School of  Mathematical Sciences, Dalian University of Technology, Dalian 116024, P.R. China}\\[1mm]
{\small E-mail: wendong@dlut.edu.cn}\\[1mm]
{\small $^\ddag$Institute of Mathematics, AMSS, Hua Loo-Keng Key
Laboratory of Mathematics,}\\[1mm]
{\small Chinese Academy of Sciences,
Beijing 100190, P.R. China}\\[1mm]
{\small E-mail: lqzhang@math.ac.cn}\\[1mm]
{\small $^\sharp$School of  Mathematical Sciences, Peking University, Beijing 100871, P.R. China}\\[1mm]
{\small E-mail: zfzhang@math.pku.edu.cn}}

\date{\today}
\maketitle

\begin{abstract}
We present some interior regularity criteria of the 3-D Navier-Stokes equations involving two components of the velocity.
These results in particular imply that if the solution is singular at one point, then at least two components of the velocity have to
blow up at the same point.
\end{abstract}

\setcounter{equation}{0}
\section{Introduction}
In this paper, we study the incompressible Navier-Stokes equations
\begin{equation}\label{eq:NS}
(NS)\left\{\begin{array}{l}
\partial_t u-\Delta u+u\cdot \nabla u+\nabla \pi=0,\\
{\rm div } u=0,
\end{array}\right.
\end{equation}
where $\big(u(x,t), \pi(x,t)\big)$ denote the velocity and the pressure of the fluid respectively.

In a seminal paper \cite{Leray}, Leray proved the global existence of weak solution with finite energy.
In two spatial dimensions, Leray weak solution is unique and regular.
In three spatial dimensions, the regularity and uniqueness of weak solution is an
outstanding open problem in the mathematical fluid mechanics. It was known that if the weak solution $u$ of (\ref{eq:NS}) satisfies
so called Ladyzhenskaya-Prodi-Serrin(LPS) type condition
\beno
\quad u\in L^q(0,T; L^p(\R^3))\quad \textrm{with} \quad \f 2 q+\f 3p\le 1, \quad p\ge 3,
\eeno
then it is regular in $\R^3\times (0,T)$, see \cite{Serrin, Giga, Struwe, ESS}, where the regularity in the class  $L^\infty(0,T; L^3(\R^3))$ was proved
by Escauriaza,  Seregin and  \v{S}ver\'{a}k \cite{ESS}.

Concerning the partial regularity of weak solution, it was started by Scheffer \cite{Sche},
and later Caffarelli, Kohn and Nirenberg \cite{CKN} showed that one dimensional Hausdorff measure of the possible singular set is zero.
The proof relies on the following small energy regularity result:
there exists some $\varepsilon_0>0$ so that if $u$ is a suitable weak solution of the Navier-Stokes equations
and satisfies
\beno
\sup_{R>0}\f 1 {R}\int_{Q_R(z)}|\na u|^2dxdt\le \varepsilon_0,
\eeno
then $u$ is regular at the point $z$ ( i.e., $u$ is bounded in a $Q_r(z)$ for some $r>0$).
Here and in what follows $z=(x,t)$, $Q_R(z)=(-R^2+t,t)\times B_R(x)$ and $B_R(x)$ is a ball of radius $r$ centered at $x$.
One could check \cite{Lin, LS, TX, GKT, Vasseur} for the simplified proof and improvements.

Recently, there are many interesting works devoted to the LPS type criterions involving the partial components of the velocity,
see \cite{CT2, CT3, CZ, KZ, PZ} and references therein.
The authors \cite{WZ2} considered the interior regularity criteria involving the partial components of the velocity.
Let
$$G(u,p,q;r)\triangleq r^{1-\frac3p-\frac2q}\|u\|_{L^q_tL^p_x(Q_r)}.$$
It was proved in \cite{WZ2} that if $(u,\pi)$ is a suitable weak solution of (\ref{eq:NS}) in $Q_1$ and satisfies
\ben\label{eq:u_3 bounded}
\sup_{0<r<1}G(u_3,p,q;r)<M \quad \textrm{for some }M>0,
\een
and 
\beno
\limsup_{r\rightarrow 0} G(u_h,p,q;r)=0,
\eeno
where $u_h=(u_1, u_2)$ and $1\le \frac3p+\frac2q<2$, $1<q\leq \infty$,
then $(0,0)$ is a regular point.

The goal of this paper is to get rid of the extra condition (\ref{eq:u_3 bounded}).
Making full use of the structure of nonlinear term and ${\rm div }u=0$,
we obtain the following interior regularity criteria involving two components of the velocity.

\begin{Theorem}\label{thm:main1}
Let $(u,\pi)$ be a suitable weak solution of (\ref{eq:NS}) in $\R^3\times (-1,0)$.
If $u$ satisfies one of the following three conditions:
\begin{itemize}
\item[1.] $u_h\in L^q_tL^p_x(Q_{r_0}), \frac3p+\frac2q=1, 2< q<\infty;$

\item[2.] $\nabla u_h\in L^q_tL^p_x(Q_{r_0}), \frac3p+\frac2q=2, 2<q<\infty;$

\item[3.] $\nabla u_h\in L^q_tL^p_x(Q_{r_0})$ and $\displaystyle\limsup_{r\rightarrow0}G(u_h,p,q;r)=0, \frac3p+\frac2q=2, 1<q\leq 2$;
\end{itemize}
for some $r_0\in (0,1)$, then $u$ is regular at $(0,0)$.
\end{Theorem}

The range of $(p,q)$ can be extended if we impose a similar condition on the velocity in a cylinder domain.
The proof relies on a new pressure decomposition formula.

\begin{Theorem}\label{thm:main2}
Let $(u,\pi)$ be a suitable weak solution of (\ref{eq:NS}) in $\R^3\times (-1,0)$. If $u$ satisfies
\beno
\limsup_{r\rightarrow0}r^{1-\frac3p-\frac2q}\Big(\int_{-r^2}^0\big(\int_{\{x;|x_h|<r,x_3\in \R\}}|u_h(x,t)|^pdx_hdx_3\big)^{\frac qp}dt\Big)^{\frac1q}=0, \quad
\eeno
where $(p,q)$ satisfies
$$1\leq\frac3p+\frac2q\le 2,\quad \frac32\leq p\leq \infty ,\quad (p,q)\neq(\infty,1),$$
then $u$ is regular at $(0,0)$.
\end{Theorem}

\begin{Remark}
An interesting consequence of Theorem \ref{thm:main1} and Theorem \ref{thm:main2} is that if the solution is singular at one point,
then at least two components of the velocity have to blow up at the same point.
\end{Remark}

\section{Suitable weak solution and $\varepsilon$-regularity criterion}

Let us first introduce the definition of suitable weak solution.

\begin{Definition} Let $\Om\subset \R^3$ and $T>0$. We say that $(u,\pi)$ is a suitable weak solution
of (\ref{eq:NS}) in $\Om_T=\Om\times (-T,0)$ if
\begin{enumerate}

\item $u\in L^{\infty}(-T,0;L^2(\Om))\cap L^2(-T,0;H^1(\Om))$ and $\pi\in L^{\frac32}(\Om_T)$;

\item the (NS) equation is satisfied in the sense of distribution;

\item the local energy inequality: for  any nonnegative $\phi\in C_c^\infty(\R^3\times\R)$
vanishing in a neighborhood of the parabolic boundary of $\Om_T$,
\beno
&&\int_{\Om}|u(x,t)|^2\phi dx+2\int_{-T}^t\int_{\Om}|\nabla u|^2\phi dxds\\
&&\leq
\int_{-T}^t\int_{\Om}|u|^2(\partial_s\phi+\triangle\phi)+u\cdot\nabla\phi(|u|^2+2\pi)dxds
\eeno
for any $t\in [-T,0]$.
\end{enumerate}
\end{Definition}

Let $(u,\pi)$ be a solution of (\ref{eq:NS}) and introduce the following scaling
\ben\label{eq:scaling}
u^{\lambda}(x, t)={\lambda}u(\lambda x,\lambda^2 t),\quad \pi^{\lambda}(x, t)={\lambda}^2\pi(\lambda x,\lambda^2 t),
\een
for any $\lambda> 0,$ then the family $(u^{\lambda}, \pi^{\lambda})$ is also a solution of (\ref{eq:NS}).
Let us introduce some invariant quantities under the scaling (\ref{eq:scaling}):
\beno
&&A(u,r,z_0)=\sup_{-r^2+t_0\leq t<t_0}r^{-1}\int_{B_r(x_0)}|u(y,t)|^2dy,\\
&&E(u,r,z_0)=r^{-1}\int_{Q_r(z_0)}|\nabla u(y,s)|^2dyds.
\eeno
We also introduce
\beno
&&G(f,p,q;r,z_0)=r^{1-\frac3p-\frac2q}\|f\|_{L^q_tL^{p}_x(Q_r(z_0))},\\
&&H(f,p,q;r,z_0)=r^{2-\frac3p-\frac2q}\|f\|_{L^q_tL^{p}_x(Q_r(z_0))},\\
&&\widetilde{G}(f,p,q;r,z_0)=r^{1-\frac3p-\frac2q}\|f-(f)_{B_r^2(x_0)}\|_{L^q_tL^{p}_x(Q_r(z_0))},\\
&&\widetilde{H}(f,p,q;r,z_0)=r^{2-\frac3p-\frac2q}\|f-(f)_{B_r^2(x_0)}\|_{L^q_tL^{p}_x(Q_r(z_0))},
\eeno
where the mixed space-time norm $\|\cdot\|_{L^q_tL^{p}_x(Q_r(z_0))}$ is defined by
\beno
\|f\|_{L^q_tL^{p}_x(Q_r(z_0))}\eqdef \Big(\int_{t_0-r^2}^{t_0}\Big(\int_{B_r(x_0)}|f(x,t)|^pdx\Big)^{\f
q p}dt\Big)^\f 1q,
\eeno
and $(f)_{B_r(x_0)}$ is the average of $f$ in the ball $B_r(x_0)$. These scaling invariant quantities will play an important role in the interior regularity theory.

For the simplicity, we denote $Q_r(0)$ by $Q_r$ and $B_r(0)$ by $B_r$, and  we will use the following notations:
$$A(u,r,(0,0))=A(u,r),\quad E(u,r,(0,0))=E(u,r).$$

Here and in what follows, we define a solution $u$ to be regular at $z_0=(x_0,t_0)$ if $u\in L^\infty(Q_r(z_0))$ for some $r>0$.
We recall the following $\varepsilon$-regularity result.


\begin{Proposition}\label{prop:small regularity-GKT}\cite{GKT}
Let $(u,\pi)$ be a suitable weak solution of (\ref{eq:NS}) in $Q_1(z_0)$ and $w=\nabla\times u$.
There exists $\varepsilon_1>0$ such that if one of the following two conditions holds,
\begin{enumerate}

\item $r^{1-\frac3p-\frac2q}\|u\|_{L^q_tL^{p}_x(Q_r(z_0))}\leq \varepsilon_1$ for any $0<r<\f12$, where $1\leq \frac3p+\frac2q\leq 2$;\vspace{0.1cm}

\item $r^{2-\frac3p-\frac2q}\|w\|_{L^q_tL^{p}_x(Q_r(z_0))}\leq \varepsilon_1$ for any $0<r<\f12$, where $2\leq \frac3{p}+\frac2{q}\leq 3$ and  $(p,q)\neq(1,\infty)$;
\end{enumerate}
then $u$ is  regular at $z_0$.
\end{Proposition}

\section{Proof of Theorem \ref{thm:main1}}

Throughout this section, we assume that $(u,\pi)$ is a suitable weak solution of (\ref{eq:NS}) in $\R^3\times (-1,0)$
and $u\in L^\infty(-1,0;L^2(\R^3))\cap L^2(-1,0;H^1(\R^3))$.

\subsection{Proof of Case 1}

In this subsection, we assume that $\nabla_h=(\partial_1,\partial_2)$ and $u_h=(u_1,u_2)\in L^q_tL^p_x(Q_{r_0})$ for some $r_0\in (0,1)$, where $\frac3p+\frac2q=1, 2< q<\infty$.
We denote by $(p',q')$ the conjugate index of $(p,q)$.

\begin{Lemma}\label{lem:local norm of nonlinear term}
It holds that for any $r\in (0,1)$,
\begin{itemize}

\item[1.] if  $\frac3l+\frac2s=\frac32, 2\leq l\leq 6$, we have
\beno
\|u\|_{L^s_tL^l_x(Q_{r})}\leq C\big(\|u\|_{L^{\infty}_tL^2_x(Q_{r})}+\|\na u\|_{L^{2}_tL^2_x(Q_{r})}\big);
\eeno

\item[2.] if  $\frac3l+\frac2s=4, 1\leq l\leq \frac32$, we have
\beno
\||u|\nabla u\|_{L^s_tL^l_x(Q_{r})}\leq C\big(\|u\|_{L^{\infty}_tL^2_x(Q_{r})}+\|\na u\|_{L^{2}_tL^2_x(Q_{r})}\big)^2.
\eeno
\end{itemize}
Here $C$ is a constant independent of $r$.
\end{Lemma}

\no{\bf Proof.}\,By scaling invariance, it suffices to consider the case of $r=1$.
By H\"{o}lder  inequality  and Sobolev interpolation inequality ( for example, see \cite{CKN}), we get
\beno
\|u\|_{L^s_tL^l_x(Q_{1})}\le C\|u\|_{L^\infty_tL^2_x(Q_{1})}^{\f {6-l} {2l}}\|u\|_{L^s_tH^1_x(Q_{1})}^{\f {3l-6} {2l}}
\le C\big(\|u\|_{L^{\infty}_tL^2_x(Q_{1})}+\|\na u\|_{L^{2}_tL^2_x(Q_{1})}\big).
\eeno
This gives the first inequality. The proof of the second inequality is similar.
\endproof


In the following, we derive the local energy inequality. We denote
\beno
G_1(f,p,q;r)\triangleq r^{3-\frac3p-\frac2q}\|f\|_{L^q_tL^{p}_x(Q_r)}.
\eeno

\begin{Lemma}\label{lem:local energy-uh}
Let $0<4r<\rho<r_0$ and $1\le p, q\le \infty$. Then we have
\beno
&&A(u,r)+E(u,r)\nonumber\\
&&\leq C\big(\frac{r}{\rho}\big)^2A(u,\rho)+C\big(\frac{\rho}{r}\big)^2G(u_h,p,q;\rho)\big(A(u,\rho)+E(u,\rho)+G_1(\nabla_h\pi,p',q';\rho)\big)\nonumber\\
&&\quad+C\big(\frac{\rho}{r}\big)^2\big(A(u,\rho)+E(u,\rho)\big)^{\frac12}\big(\tilde{H}(\pi_1,2,2;\rho)+\tilde{H}(\pi_3,2,2;\rho)\big)\\
&&\quad+C\big(\frac{\rho}{r}\big)\big(A(u,\rho)+E(u,\rho)\big)^{\frac12}G_1\big(\partial_3\pi_4,\frac{2p}{p+2},\frac{2q}{q+2};\rho\big)\nonumber,
\eeno
where the constant $C$ is independent of $r,\rho$, and $\pi_1,\pi_3$ and $\partial_3\pi_4$ is given by
\beno
&&\pi_1=\frac{1}{4\pi}\int_{\R^3}\frac{1}{|x-y|}
\sum_{i+j<6}\partial_i\partial_j\big(u_iu_j\big)dy,\quad
\pi_3=\frac{1}{2\pi}\int_{\R^3}\frac{1}{|x-y|}
\partial_3\nabla_h\cdot(-u_3u_h)dy,\\
&&\partial_3\pi_4=\frac{1}{2\pi}\int_{\R^3}\frac{1}{|x-y|}
\partial_3\partial_3\big(u_h\cdot\nabla_hu_3\big)dy.
\eeno
\end{Lemma}

\no{\bf Proof.}\,Let $\zeta$ be a cutoff function, which vanishes
outside of $Q_{\rho}$ and equals 1 in $Q_{\frac{\rho}{2}}$, and satisfies
$$|\nabla\zeta|\leq C\rho^{-1},\quad |\partial_t\zeta|+|\triangle\zeta|\leq C\rho^{-2}.$$
Define the backward heat kernel as
$$\Gamma(x,t)=\frac{1}{4\pi(r^2-t)^{3/2}}e^{-\frac{|x|^2}{4(r^2-t)}}.$$
Let $\phi=\Gamma\zeta$. Due to the local energy inequality and noting that $(\partial_t+\triangle)\Gamma=0$,
we obtain
\beno
&&\sup_t\int_{B_{\rho}}|u|^2\phi dx+\int_{Q_{\rho}}|\nabla u|^2\phi dxdt\\
&&\leq \int_{Q_{\rho}}\big(|u|^2(\triangle\phi+\partial_t\phi)-\phi u\cdot\nabla(|u|^2+2\pi-2(\pi)_{B_{\rho}})\big)dxdt\\
&&\leq \int_{Q_{\rho}}\big(|u|^2(\Gamma\triangle\zeta+\Gamma\partial_t\zeta+2\nabla\Gamma\cdot\nabla\zeta)
-\phi u\cdot\nabla(|u|^2)-\phi u\cdot\nabla(2\pi-2(\pi)_{B_{\rho}})\big)dxdt.
\eeno
It is easy to verify the following facts:
\beno
&&\Gamma(x,t)\geq C^{-1}r^{-3}\quad {\rm in} \,\, Q_r;\\
&&\phi\leq C r^{-3},\quad |\nabla\phi|\leq |\nabla\Gamma|\zeta+\Gamma|\nabla\zeta|\leq Cr^{-4};\\
&&|\Gamma\triangle\zeta|+|\Gamma\partial_t\zeta|+2|\nabla\Gamma\cdot\nabla\zeta|\leq C\rho^{-5}.
\eeno
Let
\beno
I=\int_{Q_{\rho}}\phi u\cdot\nabla(|u|^2)dxdt\triangleq I_1+I_2,
\eeno
where
\beno
I_1=\int_{Q_{\rho}}\phi u_h\cdot\nabla_h(|u|^2)dxdt,\quad I_2=\int_{Q_{\rho}}\phi u_3\cdot\nabla_3(|u|^2)dxdt.
\eeno
By H\"{o}lder inequality and Lemma \ref{lem:local norm of nonlinear term}, we have
\beno
|I_1|&\leq& Cr^{-3}\|u_h\|_{L^q_tL^p_x(Q_{\rho})}\|\nabla(|u|^2)\|_{L^{q'}_tL^{p'}_x(Q_{\rho})} \\
&\leq& Cr^{-2}\big(\frac{\rho}{r}\big)G(u_h,p,q;\rho)\big(A(u,\rho)+E(u,\rho)\big),
\eeno
and  using the facts that $\nabla\cdot u=0$ and $\frac{3}{2p'}+\frac{2}{2q'}=2$, we get by integrating by parts and H\"{o}lder inequality that
\beno
|I_2|&\leq& \big|\int_{Q_{\rho}}\phi u_3\partial_3(|u_h|^2)dxdt-2\int_{Q_{\rho}}\phi u_3^2\nabla_h\cdot u_h dxdt\big|\\
&\leq& \big|\int_{Q_{\rho}}(\phi \partial_3u_3+u_3\partial_3\phi)(|u_h|^2)dxdt-2\int_{Q_{\rho}}u_h\cdot\nabla_h( \phi u_3^2) dxdt\big|\\
&\leq& Cr^{-2}\big(\frac{\rho}{r}\big)G(u_h,p,q;\rho)\big(A(u,\rho)+E(u,\rho)+r^{-1}\|u\|_{L^{2q'}_tL^{2p'}_x(Q_{\rho})}^2\big)\\
&\leq& Cr^{-2}\big(\frac{\rho}{r}\big)^2G(u_h,p,q;\rho)\big(A(u,\rho)+E(u,\rho)\big).
\eeno
This gives that
\beno
|I|\leq Cr^{-2}\big(\frac{\rho}{r}\big)^2G(u_h,p,q;\rho)\big(A(u,\rho)+E(u,\rho)\big).
\eeno

The main trouble comes from the term including the pressure. Let
\beno
II=\int_{Q_{\rho}}\phi u\cdot\nabla(\pi-(\pi)_{B_{\rho}})dxdt\triangleq II_1+II_2,
\eeno
where
\beno
II_1=\int_{Q_{\rho}}\phi u_h\cdot\nabla_h(\pi-(\pi)_{B_{\rho}})dxdt,\quad II_2=\int_{Q_{\rho}}\phi u_3\partial_3(\pi-(\pi)_{B_{\rho}})dxdt.
\eeno
We get by H\"{o}lder inequality that
\beno
|II_1|
\leq C_0r^{-2}\big(\frac{\rho}{r}\big)G(u_h,p,q;\rho)G_1(\nabla_h\pi,p',q';\rho).
\eeno
To deal with $II_2$, recall that the pressure $\pi$ satisfies
$$-\triangle\pi=\partial_i\partial_j(u_iu_j),$$
hence,
$$
\pi=\frac{1}{4\pi}\int_{\R^3}\frac{1}{|x-y|}
\sum_{i,j=1}^3\partial_i\partial_j\big(u_iu_j\big)dy\triangleq \pi_1+\pi_2,
$$
where
\beno
\pi_1=\frac{1}{4\pi}\int_{\R^3}\frac{1}{|x-y|}
\sum_{i+j<6}\partial_i\partial_j\big(u_iu_j\big)dy,\quad \pi_2=\frac{1}{4\pi}\int_{\R^3}\frac{1}{|x-y|}
\partial_3\partial_3\big(u_3u_3\big)dy.
\eeno
We get by using $\nabla\cdot u=0$ that
\beno
\int_{Q_{\rho}}\phi u_3\partial_3\pi_2dxdt&=&\int_{Q_{\rho}}\phi u_3\partial_3\Big[\frac{1}{4\pi}\int_{\R^3}\frac{1}{|x-y|}
\partial_3\partial_3\big(u_3u_3\big)dy\Big]dxdt\\
&=&\int_{Q_{\rho}}\phi u_3\partial_3\Big[\frac{1}{2\pi}\int_{\R^3}\frac{1}{|x-y|}
\partial_3\big(u_h\cdot\nabla_hu_3-\nabla_h\cdot(u_3u_h)\big)dy\Big]dxdt\\
&=&\int_{Q_{\rho}}\phi u_3\partial_3\Big[\frac{1}{2\pi}\int_{\R^3}\frac{1}{|x-y|}
\partial_3\nabla_h\cdot(-u_3u_h)dy\Big]dxdt\\
&&+\int_{Q_{\rho}}\phi u_3\Big[\frac{1}{2\pi}\int_{\R^3}\frac{1}{|x-y|}
\partial_3\partial_3\big(u_h\cdot\nabla_hu_3\big)dy\Big]dxdt\\
&\triangleq&\int_{Q_{\rho}}\phi u_3\partial_3\pi_3dxdt+\int_{Q_{\rho}}\phi u_3\partial_3\pi_4dxdt.
\eeno
Consequently, we obtain
\beno
|II_2|&\leq& \big|\int_{Q_{\rho}}\phi u_3\partial_3(\pi_1+\pi_3)dxdt\big|+\big|\int_{Q_{\rho}}\phi u_3\partial_3\pi_4dxdt\big|\\
&\leq& C\|\partial_3(\phi u_3)\|_{L^2(Q_{\rho})}\big(\|\pi_1-(\pi_1)_{B_{\rho}}\|_{L^2(Q_{\rho})}+\|\pi_3-(\pi_3)_{B_{\rho}}\|_{L^2(Q_{\rho})}\big)\\
&&+Cr^{-3}\|u_3\|_{L^n_tL^m_x(Q_{\rho})}\|\partial_3\pi_4\|_{L^{n'}_tL^{m'}_x(Q_{\rho})},
\eeno
where  $(m',n')$ is the conjugate index of $(m,n)$ satisfying
\beno
\frac{1}{m'}=\frac1p+\frac12,\quad \frac{1}{n'}=\frac1q+\frac12,
\eeno
hence, $m=\frac{2p}{p-2}$, $n=\frac{2q}{q-2}$. Thus,
\beno
|II_2|&\leq& Cr^{-2}\big(\frac{\rho}{r}\big)^2\big(A(u,\rho)+E(u,\rho)\big)^{\frac12}\big(\tilde{H}(\pi_1,2,2;\rho)+\tilde{H}(\pi_3,2,2;\rho)\big)
\\&&+Cr^{-2}\big(\frac{\rho}{r}\big)G(u,\frac{2p}{p-2},\frac{2q}{q-2};\rho)G_1(\partial_3\pi_4,\frac{2p}{p+2},\frac{2q}{q+2};\rho).
\eeno
Noting that $\frac{3}{\frac{2p}{p-2}}+\frac{2}{\frac{2q}{q-2}}=\frac32$, we get by Lemma \ref{lem:local norm of nonlinear term} that
\beno
G(u,\frac{2p}{p-2},\frac{2q}{q-2};\rho)\leq C\big(A(u,\rho)+E(u,\rho)\big)^{\frac12},
\eeno
hence,
\beno
|II_2|&\leq& Cr^{-2}\big(\frac{\rho}{r}\big)^2\big(A(u,\rho)+E(u,\rho)\big)^{\frac12}\big(\tilde{H}(\pi_1,2,2;\rho)+\tilde{H}(\pi_3,2,2;\rho)
\big)\\
&&+Cr^{-2}\big(\frac{\rho}{r}\big)\big(A(u,\rho)+E(u,\rho)\big)^{\frac12}G_1(\partial_3\pi_4,\frac{2p}{p+2},\frac{2q}{q+2};\rho).
\eeno
Now the lemma follows by summing up the estimates of $I, II_1$ and $II_2$.\endproof\vspace{0.1cm}

The following lemma is devoted to the estimates of the pressure.

\begin{Lemma}\label{lem:pressure-uh}
Let $\pi_1,\pi_3, \partial_3\pi_4$ be as in Lemma \ref{lem:local energy-uh}. Then it holds that for $0<8r<\rho<r_0$,
\beno
&&\tilde{H}(\pi_1,2,2;r)\leq C\big(\frac{\rho}{r}\big)^\f12
G(u_h,p,q;\rho)\big(A(u,\rho)+E(u,\rho)\big)^{\frac12}+C\big(\frac{r}{\rho}\big)^{2}\tilde{H}(\pi_1,1,2;\rho),\\
&&\tilde{H}(\pi_3,2,2;r)\leq C\big(\frac{\rho}{r}\big)^\f12
G(u_h,p,q;\rho)\big(A(u,\rho)+E(u,\rho)\big)^{\frac12}+C\big(\frac{r}{\rho}\big)^{2}\tilde{H}(\pi_3,1,2;\rho),\\
&&G_1(\nabla_h\pi,p',q';r)\leq C\big(\frac{\rho}{r}\big)\big(A(u,\rho)+E(u,\rho)\big)
+C\big(\frac{r}{\rho}\big)^{\frac{3}{p'}-1}G_1(\nabla_h\pi,1,q';\rho),\\
&&G_1(\partial_3\pi_4,\frac{2p}{p+2},\frac{2q}{q+2};r)\leq C\big(\frac{\rho}{r}\big)G(u_h,p,q;\rho)E(u,\rho)^{\frac12}
+C\big(\frac{r}{\rho}\big)^{\frac{p+6}{2p}}G_1(\partial_3\pi_4,1,\frac{2q}{q+2};\rho),
\eeno
where $C$ is a constant independent of $r,\rho$.
\end{Lemma}

\no{\bf Proof.}\, Let $\zeta$ be a cut-off function,which equals 1 in $Q_{\frac{\rho}{2}}$ and vanishes outside of $Q_{\rho}$.
We decompose $\pi_1$ into $\tilde{\pi}_1+\tilde{\pi}_2$ with
$$
\tilde{\pi}_1=\frac{1}{2\pi}\int_{\R^3}\frac{1}{|x-y|}
\sum_{i+j<6}\partial_i\partial_j\big(u_iu_j\zeta^2\big).
$$
By Calderon-Zygmund inequality, we have
\beno
\int_{B_{\rho}}|\tilde{\pi}_1|^{2}dx\leq C\int_{B_{\rho}}\big(|u_h||u|\big)^{2}dx.
\eeno
Since $\tilde{\pi}_2$ is harmonic in $Q_{\frac{\rho}{2}}$, we have
\beno
\int_{B_{r}}|\tilde{\pi}_2-(\tilde{\pi}_2)_{B_{r}}|^{2}dx&\leq& C_0r^{5}\sup_{B_{\rho/4}}|\nabla\tilde{\pi}_2|^{2}\\
&\leq& C\Big(\frac{r}{\rho}\big)^{5}\rho^{-3}\big(\int_{B_{\frac{\rho}{2}}}|\tilde{\pi}_2-(\tilde{\pi}_2)_{B_{\frac{\rho}{2}}}|dx\Big)^{2}.
\eeno
Then we get by Lemma \ref{lem:local norm of nonlinear term} that
\beno
\tilde{H}(\pi_1,2,2;r)&\leq& Cr^{-\f12}\||u_h|u\|_{L^2_tL^2_x(Q_{\rho})}+C\big(\frac{r}{\rho}\big)^{2}\tilde{H}(\pi_1,1,2;\rho)\\
&\leq& Cr^{-\f12}G(u_h,p,q;\rho)\|u\|_{L^{\frac{2q}{q-2}}_tL^{\frac{2p}{p-2}}_x(Q_{\rho})}+C\big(\frac{r}{\rho}\big)^{2}\tilde{H}(\pi_1,1,2;\rho)\\
&\leq& C\big(\frac{\rho}{r}\big)^\f12
G(u_h,p,q;\rho)\big(A(u,\rho)+E(u,\rho)\big)^{\frac12}+C\big(\frac{r}{\rho}\big)^{2}\tilde{H}(\pi_1,1,2;\rho).
\eeno
The first equality of the lemma is proved. The proof of the second inequality is almost the same.
Let us turn to the proof of the third inequality.  Recall that $\pi$ satisfies
$$-\triangle\pi=\partial_i\partial_j(u_iu_j).$$
We decompose $\na_h \pi$ into $\tilde{\pi}_1+\tilde{\pi}_2$ with
$$
\tilde{\pi}_1=\frac{1}{4\pi}\int_{\R^3}\frac{1}{|x-y|}\partial_i\partial_j\big(\nabla_h(u_iu_j)\zeta^2\big)dx.
$$
By Calderon-Zygmund inequality,  we have
\beno
\int_{B_{\rho}}|\tilde{\pi}_1|^{p'}dx\leq C\int_{B_{\rho}}\big(|u||\nabla_hu|\big)^{p'}dx.
\eeno
Since $\nabla_h\tilde{\pi}_2$ is harmonic in $Q_{\frac{\rho}{2}}$, we have
\beno
\int_{B_{r}}|\tilde{\pi}_2|^{p'}dx&\leq& Cr^{3}\sup_{B_{\rho/4}}|\tilde{\pi}_2|^{p'}\\
&\leq& C\big(\frac{r}{\rho}\big)^{3}\rho^{-3p'+3}\Big(\int_{B_{\frac{\rho}{2}}}|\tilde{\pi}_2|dx\Big)^{p'}.
\eeno
Noting that $\frac{3}{p'}+\frac{2}{q'}=4$, we infer from Lemma \ref{lem:local norm of nonlinear term} that
\beno
&&G_1(\nabla_h\pi,p',q';r)\\
&&\leq G_1(\tilde{\pi}_1,p',q';r)+G_1(\tilde{\pi}_2,p',q';r)\\
&&\leq Cr^{-1}\|u\na_h u\|_{L^{q'}_tL^{p'}_x(Q_{\rho})}+C\big(\frac{r}{\rho}\big)^{\frac{3}{p'}-1}G_1(\tilde{\pi}_2,1,q';\rho)\\
&&\leq C\big(\frac{\rho}{r}\big)\big(A(u,\rho)+E(u,\rho)\big)
+C\big(\frac{r}{\rho}\big)^{\frac{3}{p'}-1}G_1(\nabla_h\pi,1,q';\rho).
\eeno
The third inequality is proved. The proof of the fourth inequality is similar. \endproof

\begin{Lemma}\label{lem:pressure estimates-fixed}
Let $\pi_1,\pi_3,\partial_3\pi_4$ be as in Lemma \ref{lem:local energy-uh}. It holds that for any $r_0\in (0,1)$,
\beno
&& \tilde{H}(\pi_1,1,2;r_0)+ \tilde{H}(\pi_3,1,2;r_0)\le C,\\
&&G_1(\nabla_h\pi,1,q';r_0)+G_1(\partial_3\pi_4,1,\frac{2q}{q+2};r_0)\le C,
\eeno
where the constant $C$ depends on $r_0$ and $\|u\|_{L^\infty(-1,0;L^2(\R^3))\cap L^2(-1,0;H^1(\R^3))}$.
\end{Lemma}

\no{\bf Proof.}\,As in the proof of Lemma \ref{lem:local norm of nonlinear term}, we have
\beno
&&\|u\|_{L^s_tL^l_x((-1,0)\times R^3)}\leq C\|u\|_{L^{\infty}_tL^2_x\cap L^{2}_tH^1_x((-1,0)\times R^3)},\quad \frac3l+\frac2s=\frac32,\quad 2\leq l\leq 6,\\
&&\|u|\nabla u|\|_{L^s_tL^l_x((-1,0)\times R^3)}\leq C\|u\|_{L^{\infty}_tL^2_x\cap L^{2}_tH^1_x((-1,0)\times R^3)}^2,\quad \frac3l+\frac2s=4,\quad 1\leq l\leq \frac32,
\eeno
from which and Calderon-Zygmund inequality, it follows that
\beno
&&\|(\pi_1,\pi_3)\|_{L^{\frac s2}_tL^{\frac l2}_x((-1,0)\times R^3)}\leq C\|u\|_{L^s_tL^l_x((-1,0)\times R^3)}^2,\quad \frac3l+\frac2s=\frac32,\quad 2< l\leq 6,\\
&&\|(\nabla\pi,\partial_3\pi_4)\|_{L^{s}_tL^{l}_x((-1,0)\times R^3)}\leq C\|u|\nabla u|\|_{L^s_tL^l_x((-1,0)\times R^3)},\quad \frac3l+\frac2s=4,\quad 1< l\leq \frac32.
\eeno
The lemma follows by taking suitable $(s,l)$ and H\"{o}lder inequality.\endproof\vspace{0.1cm}

Now we are in position to prove  Case 1 in Theorem \ref{thm:main1}. Given any $\varepsilon>0$, there exists $\rho \in (0,r_0)$ so that
\ben
G(u_h,p,q;\rho)\leq \varepsilon.
\een
Take $r$ so that $0<8r<\rho<r_0$. It follows from Lemma \ref{lem:local energy-uh} that
\beno
 &&A(u,r)+E(u,r)\\
 &&\leq C\big(\frac{r}{\rho}\big)^2A(u,\rho)+C_0\big(\frac{\rho}{r}\big)^2\varepsilon\big(A(u,\rho)+E(u,\rho)+G_1(\nabla_h\pi,p',q';\rho)\big)\\
 &&\quad+C\big(\frac{\rho}{r}\big)^2\big(A(u,\rho)+E(u,\rho)\big)^{\frac12}\big(\tilde{H}(\pi_1,2,2;\rho)+\tilde{H}(\pi_3,2,2;\rho)\big)\\
&&\quad+C\big(A(u,\rho)+E(u,\rho)\big)^{\frac12}G_1(\partial_3\pi_4,\frac{2p}{p+2},\frac{2q}{q+2};\rho)\\
 &&\leq C\big(\frac{r}{\rho}\big)^2A(u,\rho)+C\big(\frac{\rho}{r}\big)^4\big((\varepsilon+\delta)\big(A(u,\rho)+E(u,\rho)\big)+\varepsilon G_1(\nabla_h\pi,p',q';\rho)\big)\\
 &&\quad+C\delta^{-1}\big(\tilde{H}(\pi_1,2,2;\rho)^2+\tilde{H}(\pi_3,2,2;\rho)^2
 +G_1(\partial_3\pi_4,\frac{2p}{p+2},\frac{2q}{q+2};\rho)^2\big),
\eeno
where $\delta>0$ will be determined later. Let
\beno
F(r)&=&A(u;r)+E(u;r)+\varepsilon^{\frac12} G_1(\nabla_h\pi,p',q';r)\\
&&+\delta^{-\frac32}\big(\tilde{H}(\pi_1,2,2;r)^2+\tilde{H}(\pi_3,2,2;r)^2
 +G_1(\partial_3\pi_4,\frac{2p}{p+2},\frac{2q}{q+2};r)^2\big).
\eeno
Then it follows from Lemma \ref{lem:pressure-uh} that
\beno
F(r)&\leq& C\Big(\big(\frac{r}{\rho}\big)^2+(\varepsilon+\delta+\sqrt{\varepsilon})\big(\frac{\rho}{r}\big)^4+\sqrt{\delta}\Big)F(\rho) \\
&&+C\Big(\sqrt{\varepsilon}\big(\frac{\rho}{r}\big)+\big(\frac{r}{\rho}\big)^{\frac{3}{p'}-1}\Big)F(\rho) \\
&&+C\Big(\delta^{-\frac32}\big(\frac{\rho}{r}\big)^2\varepsilon^2+\big(\frac{r}{\rho}\big)^4+\big(\frac{r}{\rho}\big)^{1+\frac{6}{p}}\Big)F(\rho).
\eeno
Take $r=\theta\rho$ with $0<\theta<\f18$. The above inequality yields that
\beno
F(\theta\rho)
\leq C\big(\theta^2+\sqrt{\delta}+(\varepsilon+\delta+\sqrt{\varepsilon})\theta^{-4}
+\sqrt{\varepsilon}\theta^{-1}+\theta^{\frac{3}{p'}-1}
+\delta^{-\frac32}\theta^{-2}\varepsilon^2+\theta^{1+\frac{6}{p}}\big)F(\rho).
\eeno
We first choose $\theta$ small enough, then choose $\delta$ small, and finally choose $\varepsilon$ small enough so that
\beno
F(\theta\rho) \leq \frac12F(\rho).
\eeno
On the other hand, Lemma \ref{lem:pressure-uh} and Lemma \ref{lem:pressure estimates-fixed} imply that
\beno
F(r_0)\le C
\eeno
with $C$ depending on $r_0$ and $\|u\|_{L^\infty(-1,0;L^2(\R^3))\cap L^2(-1,0;H^1(\R^3))}$.
Then a standard iteration argument ensures that there exists $r_1>0$ such that
\beno
F(r)\leq \varepsilon_1 \quad {\rm for\,\,  any}\quad 0<r<r_1<r_0,
\eeno
which implies Case 1 of Theorem \ref{thm:main1} by Proposition \ref{prop:small regularity-GKT}.\endproof

\subsection{Proof of Case 2 and Case 3}

Let us claim that Case 2 and Case 3 in Theorem \ref{thm:main1} can be deduced from the following theorem.

\begin{Theorem}\label{thm:regularity-velocity3}
Let $(u,\pi)$ be a suitable weak solution of (\ref{eq:NS}) in $R^3\times (-1,0)$. If it satisfies
\ben\label{ass:c23}
\limsup_{r\rightarrow0}\big(H(\nabla u_h,p,q;r)+G(u_h,p,q;r)\big)=0, \quad \frac3p+\frac2q=2,\quad 1<q<\infty,
\een
then $u$ is regular at $(0,0)$.
\end{Theorem}

Indeed, the assumptions in Case 3 obviously imply (\ref{ass:c23}). Let us verify (\ref{ass:c23}) in Case 2.
In such case, $2<q<\infty$ and $\frac32<p<3$. By Poinc\'{a}re inequality, we have
$$G(u_h-(u_h)_{B_r},p,q;r)\leq CH(\nabla u_h,p,q;r)$$
for any $0<r<r_0$. Since $\frac3{p}+\frac2q=2$, we have
\beno
&&G(u_h,p,q;r)\\
&&\leq G(u_h-(u_h)_{B_{\rho}},p,q;r)+G((u_h)_{B_{\rho}},p,q;r)\\
&&\leq C\big(\frac{\rho}{r}\big)G(u_h-(u_h)_{B_{\rho}},p,q;\rho)+C\big(\frac{r}{\rho}\big)^{\frac3p-1}G(u_h,p,q;\rho)\\
&&\le C\big(\frac{\rho}{r}\big)H(\nabla u_h,p,q;\rho)+C\big(\frac{r}{\rho}\big)^{\frac3p-1}G(u_h,p,q;\rho).
\eeno
Note that $p<3$ and $\limsup_{r\rightarrow0} H(\nabla u_h,p,q;r)=0$. Then by a standard iteration,
there holds
\beno
\limsup_{r\rightarrow0} G(u_h,p,q;r)=0,
\eeno
which implies (\ref{ass:c23}). \endproof

In what follows, we assume that $\frac3p+\frac2q=2, 1<q<\infty$. We denote by $(p',q')$ the conjugate index of $(p,q)$.
To prove Theorem \ref{thm:regularity-velocity3},  we need the following  local energy inequality.

\begin{Lemma}\label{lem:local energy}
Let $0<4r<\rho<r_0$. It holds that
\beno
&&A(u,r)+E(u,r)\\
&&\leq C\big(\frac{r}{\rho}\big)^2A(u,\rho)+C\big(\frac{\rho}{r}\big)G_1(\partial_3\pi,\frac{2p}{p+1},\frac{2q}{q+1};\rho)\big(A(u,\rho)+E(u,\rho)\big)^{\frac12}\\
&&\quad+C\big(\big(\frac{\rho}{r}\big)H(\nabla u_h,p,q;\rho)+\big(\frac{\rho}{r}\big)^2G(u_h,p,q;\rho)\big)\big(A(u,\rho)+E(u,\rho)+\widetilde{H}(\pi,p',q';\rho)\big),
\eeno
where the constant $C$ is independent of $r, \rho$.
\end{Lemma}

\no{\bf Proof.}\,Since the proof is very similar to Lemma \ref{lem:local energy-uh}, we only present a sketch.
Using the same test function $\phi$ in the proof of Lemma \ref{lem:local energy-uh}, we have
\beno
&&\sup_t\int_{B_{\rho}}|u|^2\phi dx+\int_{Q_{\rho}}|\nabla u|^2\phi dxdt\\
&&\leq \int_{Q_{\rho}}\big(|u|^2(\triangle\phi+\partial_t\phi)-\phi u\cdot\nabla(|u|^2+2\pi-2(\pi)_{B_{\rho}})\big)dxdt\\
&&\leq \int_{Q_{\rho}}\big(|u|^2(\Gamma\triangle\zeta+\Gamma\partial_t\zeta+2\nabla\Gamma\cdot\nabla\zeta)
-\phi u\cdot\nabla(|u|^2)-\phi u\cdot\nabla(2\pi-2(\pi)_{B_{\rho}})\big)dxdt.
\eeno
Let
\beno
I=\int_{Q_{\rho}}\phi u\cdot\nabla(|u|^2)dxdt\triangleq I_1+I_2,
\eeno
where
\beno
I_1=\int_{Q_{\rho}}\phi u_h\cdot\nabla_h(|u|^2)dxdt,\quad I_2=\int_{Q_{\rho}}\phi u_3\cdot\pa_3(|u|^2)dxdt.
\eeno
By H\"{o}lder inequality and $\nabla\cdot u=0$, we have
\beno
|I_1|&\leq& \big|\int_{Q_{\rho}}\phi\nabla_h\cdot u_h(|u|^2)dxdt\big|+\big|\int_{Q_{\rho}}(u_h\cdot\nabla_h\phi) (|u|^2)dxdt\big|\\
&\leq& Cr^{-2}\big(\big(\frac{\rho}{r}\big)H(\nabla u_h,p,q;\rho)+\big(\frac{\rho}{r}\big)^2G(u_h,p,q;\rho)\big)G(u,2p',2q';\rho)^2,
\eeno
and noting that $\pa_3|u|^2\le |\na u_h||u|$, we get
\beno
|I_2|\leq Cr^{-2}\big(\frac{\rho}{r}\big)H(\nabla u_h,p,q;\rho)G(u,2p',2q';\rho)^2,
\eeno
which along with Lemma \ref{lem:local norm of nonlinear term} imply that
\beno
|I| \leq Cr^{-2}\big(\big(\frac{\rho}{r}\big)H(\nabla u_h,p,q;\rho)+\big(\frac{\rho}{r}\big)^2G(u_h,p,q;\rho)\big)\big(A(u,\rho)+E(u,\rho)\big).
\eeno
Let
\beno
II=\int_{Q_{\rho}}\phi u\cdot\nabla(\pi-(\pi)_{B_{\rho}})dxdt\triangleq II_1+II_2,
\eeno
where
\beno
II_1=\int_{Q_{\rho}}\phi u_h\cdot\nabla_h(\pi-(\pi)_{B_{\rho}})dxdt,\quad II_2=\int_{Q_{\rho}}\phi u_3\partial_3(\pi-(\pi)_{B_{\rho}})dxdt.
\eeno
We have by H\"{o}lder inequality and Lemma \ref{lem:local norm of nonlinear term} that
\beno
|II_1|&\leq& \big|\int_{Q_{\rho}}\phi\nabla_h\cdot u_h((\pi-(\pi)_{B_{\rho}})dxdt\big|+\big|\int_{Q_{\rho}}(u_h\cdot\nabla_h\phi) (\pi-(\pi)_{B_{\rho}})dxdt\big|\\
&\leq& Cr^{-2}\big(\big(\frac{\rho}{r}\big)H(\nabla u_h,p,q;\rho)+\big(\frac{\rho}{r}\big)^2G(u_h,p,q;\rho)\big)\widetilde{H}(\pi,p',q';\rho),\\
|II_2|&\leq& Cr^{-2}\big(\frac{\rho}{r}\big)G_1(\partial_3\pi,\frac{2p}{p+1},\frac{2q}{q+1};\rho)G(u,2p',2q';\rho)\\
&\le& Cr^{-2}\big(\frac{\rho}{r}\big)G_1(\partial_3\pi,\frac{2p}{p+1},\frac{2q}{q+1};\rho)\big(A(u,\rho)+E(u,\rho)\big)^{\frac12}.
\eeno
The lemma follows by summing up the estimates of $I, II_1$ and $II_2$.\endproof

The proof of the following lemma is similar to Lemma \ref{lem:pressure-uh}. So, we omit the details.

\begin{Lemma}\label{lem:pressure-nabla uh}
It holds that for any $0<8r<\rho<r_0$,
\beno
&&\tilde{H}(\pi,p',q';r)\leq C\big(\frac{\rho}{r}\big)
\tilde{G}(u,2p',2q';\rho)^2+C\big(\frac{r}{\rho}\big)^{\frac{3}{p'}}\tilde{H}(\pi,1,q';\rho),\\
&&G_1(\partial_3\pi,\frac{2p}{p+1},\frac{2q}{q+1};r)
\leq C\big(\frac{\rho}{r}\big)^{\frac12}\tilde{G}(u,2p',2q';\rho)H(\nabla u_h,p,q;\rho)\\
&&\qquad\qquad\qquad+C\big(\frac{r}{\rho}\big)^{1+\frac{3}{2p}}G_1(\partial_3\pi,1,\frac{2q}{q+1};\rho),
\eeno
where the constant $C$ is independent of $r, \rho$.
\end{Lemma}

Now let us turn to prove Theorem \ref{thm:regularity-velocity3}. By the assumption, given any $\varepsilon>0$,
there exists $\rho \in (0,r_0)$ so that
\beno
H(\nabla u_h,p,q;\rho)+G(u_h,p,q;\rho)\leq \varepsilon.
\eeno
Take $r>0$ so that $0<8r<\rho<r_0$. It follows from Lemma \ref{lem:local energy} that
\beno
A(u,r)+E(u,r)&\leq& C\big(\frac{r}{\rho}\big)^2A(u,\rho)+C\delta^{-1}G_1(\partial_3\pi,\frac{2p}{p+1},\frac{2q}{q+1};\rho)^2\\
 &&+C\big(\frac{\rho}{r}\big)^2\big((\varepsilon+\delta)(A(u,\rho)+E(u,\rho))+\varepsilon\widetilde{H}(\pi,p',q';\rho)\big),
\eeno
where $\delta>0$  will be determined later. Let
\beno
F(r)=A(u,b;r)+E(u,b;r)+\varepsilon^{1/2} \tilde{H}(\pi,p',q';r)+\delta^{-\frac32}G_1(\partial_3\pi,\frac{2p}{p+1},\frac{2q}{q+1};r)^2.
\eeno
Then it follows from Lemma \ref{lem:pressure-nabla uh} that
\beno
F(r)&\leq& C\Big(\big(\frac{r}{\rho}\big)^2+\sqrt{\delta}+(\varepsilon+\delta+\sqrt{\varepsilon})\big(\frac{\rho}{r}\big)^2\Big)F(\rho) \\
&&+C\Big(\sqrt{\varepsilon}\big(\frac{\rho}{r}\big)+\big(\frac{r}{\rho}\big)^{\frac{3}{p'}}\Big)F(\rho)
+C\Big(\delta^{-\frac32}\big(\frac{\rho}{r}\big)\varepsilon^2+\big(\frac{r}{\rho}\big)^{2+\frac{3}{p}}\Big)F(\rho).
\eeno
Take $r=\theta\rho$ with $0<\theta<\f18$. The above inequality yields that
\beno
F(\theta\rho)\leq C\big(\theta^2+\sqrt{\delta}+(\varepsilon+\delta+\sqrt{\varepsilon})\theta^{-2}
+\sqrt{\varepsilon}\theta^{-1}+\theta^{\frac{3}{p'}}
+\delta^{-\frac32}\theta^{-1}\varepsilon^2+\theta^{2+\frac{3}{p}}\big)F(\rho).
\eeno
We first choose $\theta$ small enough, then choose $\delta$ small, finally choose $\varepsilon$ small enough so that
\beno
F(\theta\rho)\leq \frac12F(\rho).
\eeno
On the other hand, it is easy to see that
\beno
F(r_0)\le C
\eeno
with $C$ depending on $r_0$ and $\|u\|_{L^\infty(-1,0;L^2(\R^3))\cap L^2(-1,0;H^1(\R^3))}$. Then a standard iteration argument ensures that there exists $r_1>0$ such that
\beno
F(r)\leq \varepsilon_1 \quad {\rm for\,\, all}\quad 0<r<r_1<r_0.
\eeno
which implies Theorem \ref{thm:regularity-velocity3} by Proposition \ref{prop:small regularity-GKT}.\endproof

\section{Proof of Theorem \ref{thm:main2}}

Throughout this section, we assume that $(u,\pi)$ be a suitable weak solution of (\ref{eq:NS}) in $\R^3\times (-1,0)$.
Let us first introduce some notations.

Let $B_r^2=\big\{(x_1,x_2);|(x_1,x_2)|\leq r\big\}$, $B_r^*=B_r^2\times R=\big\{(x_1,x_2,x_3);(x_1,x_2)\in B_r^2, x_3\in R\big\}$,
and $Q_r^*=B_r^*\times (-r^2,0)$. Moreover, $Q_r^*(z_0)=(-r^2+t_0,t_0)\times B_r^*(x_0)$, $B_r^*(x_0)=B_r^2(x_0)\times R$
and $B_r^2(x_0)$ is a ball of radius $r$ centered at the horizontal part of $x_0$. For the simplicity, we
denote $Q_r^*(0)$ by $Q_r^*$ and $B_r^*(0)$ by $B_r^*$.
As in Section 2, we will still use the notations like $A(u,r),  E(u,r),$  $G(f,p,q;r),$ $H(f,p,q;r), $ $\widetilde{G}(f,p,q;r),$ $\widetilde{H}(\pi,p,q;r)$ etc.
The differences are that here the integral domain is replaced by $Q_r^*$ or $B_r^*$, and the mean value in $\widetilde{G}, \widetilde{H}$ is taken only on $B_r^2.$
We denote by $(p',q')$ the conjugate index of $(p,q)$.

\begin{Lemma}\label{prop:local}
Let $0<4r<\rho<r_0$ and $1\le p, q\le \infty$. We have
\beno
&&A(u;r)+E(u;r)\nonumber\\
&&\leq
C\big(\frac{r}{\rho}\big)A(u;\rho)+C\big(\frac{\rho}{r}\big)^2G(u_h,p,q;\rho)
\Big(G(u,2p',2q';\rho)^2+\tilde{H}(\pi,p',q';\rho)\Big),
 \eeno
where $C$ is a constant independent of $r,\rho$.
\end{Lemma}

\no{\bf Proof.}\,Let $\zeta$ be a cutoff function, which vanishes
outside of $Q^*_{\rho}$ and equals 1 in $Q^*_{\frac{\rho}{2}}$, and satisfies
$$|\nabla\zeta|\leq C_0\rho^{-1},\quad |\partial_t\zeta|+|\triangle\zeta|\leq C_0\rho^{-2}.$$
Define the backward heat kernel as
$$\Gamma(x,t)=\frac{1}{4\pi(r^2-t)}e^{-\frac{|x_h|^2}{4(r^2-t)}}.$$
Taking the test function $\phi=\Gamma\zeta$ in the local energy inequality,
and noting $(\partial_t+\triangle_h)\Gamma=0$, where $\triangle_h=\partial_{x_1}^2+\partial_{x_2}^2$,
we obtain
\beno
&&\sup_t\int_{B^*_{\rho}}|u|^2\phi dx+\int_{Q^*_{\rho}}|\nabla u|^2\phi dxdt\\
&&\leq \int_{Q^*_{\rho}}\big[|u|^2(\triangle\phi+\partial_t\phi)+u\cdot\nabla\phi(|u|^2+2\pi-2(\pi)_{B^2_{\rho}})\big]dxdt\\
&&\leq \int_{Q^*_{\rho}}\big[|u|^2(\Gamma\triangle\zeta+\Gamma\partial_t\zeta+2\nabla\Gamma\cdot\nabla\zeta)
+|\nabla\phi||u_h|(|u|^2+2|\pi-(\pi)_{B^2_{\rho}}|)\big]dxdt.
\eeno
It is easy to verify that
\beno
&&\Gamma(x,t)\geq C_0^{-1}r^{-2}\quad {\rm in} \,\, Q_r^*,\\
&&|\nabla\phi|\leq |\nabla\Gamma|\zeta+\Gamma|\nabla\zeta|\leq C_0r^{-3},\\
&&|\Gamma\triangle\zeta|+|\Gamma\partial_t\zeta|+2|\nabla\Gamma\cdot\nabla\zeta|\leq C_0\rho^{-4},
\eeno
from which and H\"{o}lder inequality, it follows that
\beno
 &&A(u,r)+E(u,r)\\
 &&\leq C\big(\frac{r}{\rho}\big)A(u,\rho)+C\big(\frac{\rho}{r}\big)^2\rho^{-2}\int_{Q^*_{\rho}}(|u_h||u|^2
 +|u_h||\pi-(\pi)_{B^2_{\rho}}|)dxdt\\
 &&\leq C\big(\frac{r}{\rho}\big)A(u,\rho)+C\big(\frac{\rho}{r}\big)^2G(u_h,p,q;\rho)
\big(G(u,2p',2q';\rho)^2+\tilde{H}(\pi,p',q';\rho)\big).
\eeno
This completes the proof of the lemma. \endproof

In the sequel, we assume that $(p,q)$ satisfies
\beno
\frac3p+\frac2q=2,\quad \frac32\leq p<\infty.
\eeno

%
%

\begin{Lemma}\label{lem:magnetic}
For any $0<r<r_0$, we have
\beno
G(u,2p',2q';r)^2\leq C\big(E(u,2r)+A(u,2r)\big),
\eeno
where $C$ is a constant independent of $r$.
\end{Lemma}

\no{\bf Proof.}\, Recall a well-known Sobolev's interpolation inequality (for example, see \cite{CKN}):
\ben\label{eq:Sobolev}
\int_{\R^3}|f|^{\ell}\leq C\Big(\int_{\R^3}|\nabla f|^2dx\Big)^a
\Big(\int_{\R^3}|f|^2dx\Big)^{\frac \ell 2-a},
\een
where $2\leq \ell\leq 6$ and $a=\frac34(\ell-2)$. Applying (\ref{eq:Sobolev}) with $\ell=2p'$ (Note that $2p'\leq 6$ since $p\geq \f 32$) and a suitable localization,
we get
\beno
&&G(u,2p',2q';r)^2=r^{\f 3p+\f 2q-3}\|u\|^2_{L^{2p',2q'}(Q_r^*)}\nonumber\\
&&\leq  Cr^{-1}\Big\{\int_{-r^2}^0\big[\big(\int_{B_{2r}^*}|\nabla u|^2\big)^a(\int_{B_{2r}^*}|u|^2\big)^{p'-a}
+r^{-2a}\big(\int_{B_{2r}^*}|u|^2\big)^{p'}\big]^{\frac{q'}{p'}}dt\Big\}^{\frac{1}{q'}}\nonumber\\
&&\leq Cr^{-1}\Big\{\int_{-r^2}^0\big(\int_{B_{2r}^*}|\nabla u|^2\big)^{\frac{aq'}{p'}}
\big(\int_{B_{2r}^*}|u|^2\big)^{q'(1-\frac{a}{p'})}dt+
r^{-\frac{2aq'}{p'}+2}\big(\sup_{t}\int_{B_{2r}^*}|u|^2\big)^{q'}\Big\}^{\frac{1}{q'}},
\eeno
then the lemma follows by noting that $\f {aq'} {p'}=1$ and
$-\frac{2a}{p'}+\f 2 {q'}=-\f 3 p-\f 2 q+2=0$.
\endproof\vspace{0.1cm}

In the following, we will introduce a new pressure decomposition formula in a cylinder domain
based on the following properties of harmonic function, which is new even for harmonic function to our knowledge.

\begin{Lemma}\label{lem:harmonic function}
Let $f$ be a harmonic function in a cubic $D_1\subset \R^3$. Let
\beno
&&P_3 f(x_h)=\frac12\int_{-1}^1f(x_1,x_2,x_3)dx_3,\\
&&P_h f(x_3)=\frac14\int_{-1}^1\int_{-1}^1f(x_1,x_2,x_3)dx_h.
\eeno
Then it holds that
\beno
&&\sup_{x\in B_{\frac12}}|\nabla_3 f|\leq C\int_{B_1}|f(x)-P_3 f(x_h)|dx,\label{eq:harmonic property1}\\
&&\sup_{x\in B_{\frac12}}|\nabla_h f|\leq C\int_{B_1}|f(x)-P_h f(x_3)|dx.\label{eq:harmonic property2}
\eeno
\end{Lemma}

\no{\textbf{Proof.}\,
For  $|h|\leq \frac15,$ let
$$f^h(x)=f^h(x_1,x_2,x_3)=f(x_1,x_2,x_3+h),$$
and $g(x)=f(x)-P_3 f(x_1,x_2)$. It is easy to see that
\beno
f(x)-f^h(x)=g(x)-g^h(x),\quad  x\in B_{\frac45}.
\eeno
Since  $f$ is a harmonic function in $D_1$, we have
\beno
\triangle\big(f(x)-f^h(x)\big)=\triangle\big(g(x)-g^h(x)\big)=0,\quad  x\in B_{\frac45}.
\eeno
The gradient estimate of harmonic function yields that
\beno
&&\sup_{B_{\frac12}}|\partial_3 f-\partial_3 f^h|\leq C \sup_{B_{\frac34}}|f- f^h|\leq C \sup_{B_{\frac34}}|g- g^h|,\\
&&\sup_{B_{\frac34}}|g- g^h|\leq C\int_{B_{\frac45}}|g-g^h|dx\leq C\int_{B_1}|f-P_3 f|dx.
\eeno
This proves that for any $|h|\leq \frac15,$
\ben\label{eq:estimate f harmonic function}
\sup_{B_{\frac12}}|\partial_3 f-\partial_3 f^h|\leq C\sup_{B_{\frac34}}|f- f^h|\leq C\int_{B_1}|f-P_3 f|dx.
\een
The second inequality of (\ref{eq:estimate f harmonic function}) implies by Mean value theorem that given $x\in B_{1/2}$,
there exists $h=h(x)$ with $|h|\le \f15$ so that
\beno
|\partial_3f(x_1,x_2,x_3+h)|\leq C\int_{B_1}|f-P_3 f|dx,
\eeno
which along with (\ref{eq:estimate f harmonic function}) gives the first inequality of the lemma.
The proof of the second inequality of the lemma is similar.\endproof

Let
\beno
\tilde{H}(\pi,p',q';r)=r^{2-\frac3{p'}-\frac2{q'}}\Big(\int_{-r^2}^0\big(\int_{B_r^*}|\pi-P_{h,r}\pi(x_3)|^{p'}dx_hdx_3\big)^{\frac{q'}{p'}}dt\Big)^{1/{q'}},
\eeno
where
\beno
P_{h,r} \pi(x_3)=\frac1{|B_r^2|}\int_{B_r^2}\pi(x_h,x_3)dx_h.
\eeno

\begin{Lemma}\label{lem:pressure1}
For any $0<8r<\rho<r_0$, it holds that
\beno
\tilde{H}(\pi,p',q';r)\leq C\big(\frac{\rho}{r}\big)
\tilde{G}(u,2p',2q';\rho)^2+C\big(\frac{r}{\rho}\big)^{\frac{2}{p'}}\tilde{H}(\pi,p',q';\rho),
\eeno
where $C$ is a constant independent of $r, \rho$.
\end{Lemma}

\no{\bf Proof.}\,Recall that the pressure $\pi$ satisfies
\beno
\triangle \pi=-\partial_i\partial_j(u_iu_j).
\eeno
Let $ \pi=\pi_1+\pi_2$ where $\pi_1$ is defined by
\beno
\triangle \pi_1=-\partial_i\partial_j(u_iu_j\chi(x_h)),
\eeno
here $\chi(x_h)$ is a smooth function with $\chi(x_h)=1$ for $|x_h|\le \f \rho 2$ and $\chi(x_h)=0$ for $|x_h|\ge \rho$.
So, $\pi_2$ is harmonic in $B_{\frac{\rho}{2}}^*$.

Due to $p'>1$, by Calderon-Zygmund inequality  we have
\ben\label{eq:pi-1}
\tilde{H}(\pi_1,p',q';\rho)\leq C\tilde{G}(u,2p',2q';\rho)^2.
\een
We denote $B_{\rho,k}^*=B_{\rho}^2\times(-k\rho,k\rho).$
Since $\pi_2$ is harmonic in $B_{\frac{\rho}{2}}^*$, we have
\beno
\int_{B_{r,k}^*}|\pi_2-(\pi_2)_{B_r^2}|^{p'}dx \leq Cr^{3+p'}\sum_{j=1}^{2k}\sup_{B_{\sqrt{2}r}(z_j)}|\nabla_h\pi_2|^{p'},
\eeno
where $z_j\in \big\{x; x_h=0,x_3=\ell r+\frac12 r-kr,0\leq \ell<2k\big\}$.
We infer from Lemma \ref{lem:harmonic function} that
\beno
\sup_{B_{\sqrt{2}r}(z_j)}|\nabla_h\pi_2|^{p'}\leq C\rho^{-3-p'}\int_{B_{\frac{\rho}{2}}(z_j)}|\pi_2-P_{h,\frac{\rho}{2}}\pi_2|^{p'}dx.
\eeno
Note that the ball $B_{\frac{\rho}{2}}(z_j)$ intersects each other at most $C\f \rho r$ times. We infer that
\beno
\int_{B_{r,k}^*}|\pi_2-(\pi_2)_{B_r^2}|^{p'}dx&\leq& C\big(\frac{r}{\rho}\big)^{2+p'}\int_{B_{\frac{\rho}{2}}^*}|\pi_2-P_{h,\frac{\rho}{2}}\pi_2|^{p'}dx,
\eeno
and letting $k\rightarrow\infty$, we get
\beno
\tilde{H}(\pi_2,p',q';r)
\leq C\big(\frac{r}{\rho}\big)^{3-\frac{1}{p'}-\frac{2}{q'}}\tilde{H}(\pi_2,p',q';\rho).
\eeno
which along with (\ref{eq:pi-1}) gives
\beno
\tilde{H}(\pi,p',q';r)
&\leq& \tilde{H}(\pi_1,p',q';r)+\tilde{H}(\pi_2,p',q';r)\\
&\leq& C\big(\frac{r}{\rho}\big)^{2-\frac{3}{p'}-\frac{2}{q'}}\tilde{G}(u,b,2p',2q';\rho)^2
+C\big(\frac{r}{\rho}\big)^{3-\frac{1}{p'}-\frac{2}{q'}}\tilde{H}(\pi,p',q';\rho).
\eeno
The proof is finished.\endproof

Now we are in position to prove Theorem \ref{thm:main2}. Let
\beno
F(r)=A(u,r)+E(u,r)+\varepsilon^{\frac12} \tilde{H}(\pi,p',q';r).
\eeno
Take $(r,\rho,\kappa)$ so that $0<8r<\rho$ and $8\rho<\kappa<r_0$ and
\beno
G(u_h, p, q; \rho)\le \varepsilon.
\eeno
It suffices to consider the case
\beno
\frac3p+\frac2q=2,\quad \frac32\leq p\leq \infty ,\quad (p,q)\neq(\infty,1),
\eeno
since the other cases can be reduced to such case by H\"{o}lder inequality.

We know from Lemma \ref{lem:magnetic} and Lemma \ref{lem:pressure1} that
\beno
\tilde{H}(\pi,p',q';\rho)&\leq& C\big(\frac{\kappa}{\rho}\big)
\tilde{G}(u,b,2p',2q';\kappa)^2+C\big(\frac{\rho}{\kappa}\big)^{\frac{2}{p'}}\tilde{H}(\pi,p',q';\kappa)\\
&\leq& C\big(\frac{\kappa}{\rho}\big)
F(\kappa)+C\big(\frac{\rho}{\kappa}\big)^{\frac{2}{p'}}\tilde{H}(\pi,p',q';\kappa),
\eeno
which along with Lemma \ref{prop:local} and Lemma \ref{lem:pressure1} gives
\beno
F(r)&\leq& C\big(\frac{r}{\rho}\big)A(u,\rho)+C\big(\frac{\rho}{r}\big)^2\big(\varepsilon G(b,2p',2q';\rho)^2+\varepsilon\tilde{H}(\pi,p',q';\rho)\big)\\
&&+\varepsilon^{\frac12}\tilde{H}(\pi,p',q';r)\\
&\leq& C\big(\frac{r}{\rho}\big)F(\rho)+C\varepsilon\big(\frac{\rho}{r}\big)^2F(\rho)
+C\big(\frac{\rho}{r}\big)^2\varepsilon\tilde{H}(\pi,p',q';\rho)+\varepsilon^{\frac12}\tilde{H}(\pi,p',q';r)\\
&\leq& C\big(\frac{r}{\rho}\big)F(\rho)+C_0\varepsilon\big(\frac{\rho}{r}\big)^2 F(\rho)
+C\big(\frac{\rho}{r}\big)^2\varepsilon\big(\frac{\kappa}{\rho}\big)F(\kappa)
+C\big(\frac{\rho}{r}\big)^2\varepsilon\big(\frac{\rho}{\kappa}\big)^{\frac{2}{p'}}\tilde{H}(\pi,p',q';\kappa)\\
&&+C\varepsilon^{\frac12}\big(\frac{\kappa}{r}\big)F(\kappa)+C\varepsilon^{\frac12}
\big(\frac{r}{\kappa}\big)^{\frac{2}{p'}}\tilde{H}(\pi,p',q';\kappa)\\
&\leq& C\big((\frac{r}{\rho})(\frac{\kappa}{\rho})+(\frac{\rho}{r})^2\varepsilon(\frac{\kappa}{\rho})
+(\frac{\rho}{r})^2\varepsilon^{\frac12}\big(\frac{\rho}{\kappa}\big)^{\frac{2}{p'}}
+\varepsilon^{\frac12}(\frac{\kappa}{r})+\big(\frac{r}{\kappa}\big)^{\frac{2}{p'}}\big)F(\kappa).
\eeno
Take $r=\theta^2\rho,\,\rho=\theta\kappa$ with $0<\theta<\f18$. The above inequality yields that
\beno
F(r)\leq C\big(\theta+\varepsilon\theta^{-5}+\varepsilon^{\frac12}\theta^{-4+\frac{2}{p'}}
+\varepsilon^{\frac12}\theta^{-3}+\theta^{\frac{6}{p'}}\big)F(\kappa).
\eeno
Choose $\theta$ small enough, and then choose $\varepsilon$ small enough so that
\beno
C\big(\theta+\varepsilon\theta^{-5}+\varepsilon^{\frac12}\theta^{-4+\frac{2}{p'}}
+\varepsilon^{\frac12}\theta^{-3}+\theta^{\frac{6}{p'}}\big)\leq \frac12.
\eeno
This gives the following iterative inequality
$$F(\theta^{2}\kappa)\leq \frac12 F(\kappa).$$
On the other hand, it is easy to see that
\beno
F(R)\le C
\eeno
with $C$ depending on $R$ and $\|u\|_{L^\infty(-1,0;L^2(\R^3))\cap L^2(-1,0;H^1(\R^3))}$. Indeed, since $\pi$ satisfies
\beno
-\triangle\pi=\partial_i\partial_j({u}_i{u}_j),
\eeno
by Calderon-Zygmund inequality and interpolation inequality, we get
\beno
\|\pi\|_{L^{q'}(-T,0;L^{p'}(\R^3))}&\le& C_0\|u\|_{L^{2q'}(-T,0;L^{2p'}(\R^3))}^2\\
&\le& C_0 \|u\|_{L^\infty(-T,0;L^2(\R^3))\cap L^2(-T,0;H^1(\R^3))}^2.
\eeno
Then a standard iteration argument ensures that there exists $r_1>0$ such that
\beno
F(r)\leq \varepsilon \quad {\rm for\,\, all}\quad 0<r<r_1.
\eeno
which implies Theorem \ref{thm:main2} by Proposition \ref{prop:small regularity-GKT}.\endproof

\bigskip

\noindent {\bf Acknowledgments.} W. Wang was supported by NSFC
11301048. L. Zhang was partially supported by the innovation program
at CAS and National Basic Research Program of China under grant
2011CB808002.  Z. Zhang was partially supported by NSF of China
under Grant 11071007, Program for New Century Excellent Talents.


\begin{thebibliography}{WWW}


\bibitem{CKN}  L. Caffarelli, R. Kohn and L. Nirenberg, {\it Partial regularity of suitable weak solutions of the Navier-Stokes equations}, Comm. Pure Appl. Math., 35(1982), 771-831.

\bibitem{CT2}  C. Cao and E. S. Titi, {\it Regularity criteria for the three dimensional Navier-Stokes equations},
Indiana Univ. Math. J.,  57 (2008),  2643-2662.

\bibitem{CT3}  C. Cao and E. S. Titi,  {\it Global regularity criterion of the 3D Navier-Stokes equations involving one entry of the velocity  gradient tensor},
Arch. Ration. Mech. Anal.,  202(2011), 919-932.

\bibitem{CZ} J. Y. Chemin and P. Zhang, {\it On the critical one component regularity for 3-D Navier-Stokes system},  arXiv:1310.6442.

\bibitem{ESS} L. Escauriaza, G. A. Seregin and V. \v{S}ver\'{a}k,
{\it $L^{3,\infty}$ solutions to the Navier-Stokes equations and backward uniqueness},
 Russ. Math. Surveys, 58(2003), 211-250.

\bibitem{Giga} Y. Giga,  {\it Solutions for semilinear parabolic equations in $L^p$ and regularity of weak solutions of the Navier-Stokes system}, J. Differential Equations, 62(1986), 186-212.

\bibitem{GKT} S. Gustafson, K. Kang and T.-P. Tsai,  {\it Interior regularity criteria for suitable weak solutions of the Navier-Stokes equations,} Comm. Math. Phys., 273(2007), 161-176.

\bibitem{Ku} I. Kukavica, {\it On partial regularity for the Navier-Stokes equations}, Discrete Contin. Dyn. Syst., 21 (2008), 717-728.

\bibitem{KZ} I. Kukavica and M. Ziane, {\it  One component regularity for the Navier-Stokes equation},  Nonlinearity,  19 (2006),  453-469.

\bibitem{LS} O. A. Ladyzhenskaya and G. A. Seregin,
{\it On partial regularity of suitable weak solutions to the three-dimensional Navier-Stokes equations}, J. Math. Fluid Mech., 1(1999), 356-387.

\bibitem{Leray} J. Leray, {\it Sur le mouvement d'un liquids visqeux emplissant
              l'espace}, Acta Math., 63(1934), 193-248.

\bibitem{Lin} F. H. Lin,  {\it A new proof of the Caffarelli-Kohn-Nirenberg theorem}, Comm. Pure Appl. Math., 51(1998), 241-257.

\bibitem{NRS} J. Ne\v{c}as, M. R\.{u}\v{z}i\v{c}ka and V. \v{S}ver\'{a}k, {\it On Leray's self-similar solutions of the Navier-Stokes equations},
 Acta Math., 176(1996), 283-294.

\bibitem{PZ} M. Pokorn\'{y} and Y. Zhou, {\it On the regularity of the solutions of the Navier-Stokes equations via one velocity component},
Nonlinearity,  23 (2010), 1097-1107.

\bibitem{Sche} V. Scheffer, {\it Partial regularity of solutions to the Navier-Stokes equations,} Pacific J. Math., 66 (1976), 535-552.


\bibitem{Serrin} J. Serrin, {\it The initial value problem for the Navier-stokes
                equations}, in Nonlinear problems(R. E. Langer Ed.),
                pp.69-98, Univ. of Wisconsin Press, Madison, 1963.

\bibitem{Struwe} M. Struwe, {\it On partial regularity results for the Navier-Stokes equations,} Comm. Pure Appl. Math., 41(1988), 437-458.

\bibitem{TX} G. Tian and Z. Xin, {\it Gradient estimation on Navier-Stokes equations}, Comm. Anal. Geom., 7(1999), 221-257.

\bibitem{Vasseur} A. Vasseur, {\it A new proof of partial regularity of solutions to Navier-Stokes equations},
Nonlinear Differential Equations Appl., 14(2007), 753-785.


\bibitem{WZ2} W. Wang and Z. Zhang, {\it  On the interior regularity criterion and the number of singular points to the  Navier-Stokes equations},
 J. Anal. Math., 123(2014), 139--170.



\end{thebibliography}
\end{document}